\documentclass[11pt]{article}
\usepackage{amssymb}
\usepackage{latexsym}

\newcommand{\reals}{\mathbb R}

\begin{document}

\title {\bf   Etude  des $(n+1)$-tissus de courbes  en dimension $n$ }

\author{   Jean-Paul Dufour  et Daniel Lehmann\\}
\maketitle

\begin{abstract} 
For $(n+1)$-webs by curves in an ambiant $n$-dimensional manifold, we first define a generalization of the well known Blaschke curvature of the dimension two, which vanishes iff the web has the maximum possible rank which is one.  But, contrary to the   dimension two where all 3-webs of rank one are locally isomorphic, we prove that there are infinitely many classes of  isomorphism for germs of  4-webs by curves of rank one in the  dimension three : 
 we  provide a procedure for building all of them, up to isomorphism,  and give examples of invariants  of these classes allowing  in particular to distinguish the so-called quadrilateral webs among them. 
\end{abstract}



\section{Introduction}

 Le contexte est  holomorphe ou analytique r\'eel.

Dans une vari\'et\'e $\cal U$ de dimension $n$, on va s'int\'eresser aux $(n+1)$-tissus de courbes (codimension $(n-1)$),  localement d\'efinis  par $(n+1)$ feuilletages ${\cal F}_\lambda$, $(1\leq \lambda \leq n+1)$, dont on note $T( {\cal F}_\lambda)$ le fibr\'e tangent et $N( {\cal F}_\lambda)$ le fibr\'e normal.
$$0\to T( {\cal F}_\lambda)\to T{\cal U} \buildrel{\pi_\lambda}\over{\longrightarrow} N( {\cal F}_\lambda)\to 0.$$
On supposera toujours 
les feuilletages du tissu en \emph{position g\'en\'erale} : notant $V_\lambda$ un champ de vecteurs engendrant $T( {\cal F}_\lambda)$, $n$ quelconques parmi les $n+1$ champs de vecteurs $V_\lambda$ seront toujours suppos\'es  lin\'eairement ind\'ependants.

Rappelons  

- qu'une $p$-forme ${\cal F}$-basique pour un feuilletage ${\cal F}$ de codimension $\geq p$ est une $p$-forme $\eta$ sur la vari\'et\'e ambiante telle que $\iota_v\eta=0$ et $L_v\eta=0$ quel que soit le champ de vecteurs $v$ tangent au feuilletage, $\iota_v$ et $L_v(=\iota_v\circ d+d\circ \iota_v)$ d\'esignant respectivement le produit int\'erieur et la d\'eriv\'ee de Lie, 

- et qu'une $(n-1)$-relation ab\'elienne (ou \emph{relation ab\'elienne} tout court en abr\'eg\'e) au voisinage d'un point de $V$ est la donn\'ee d'une famille $(\eta_\lambda)_{1\leq \lambda\leq n+1}$ de $(n-1$)-formes, telle que $\sum_\lambda\eta_\lambda=0$, chaque $\eta_\lambda$ \'etant ${\cal F}_\lambda$-basique.

 On sait que le rang  d'un tel $(n+1)$-tissu en courbes (c'est-\`a-dire la dimension maximum de l'espace des germes de $(n-1)$-relations ab\'eliennes en un point) est 0 ou 1, en vertu   d'une formule tr\`es g\'en\'erale de Damiano ([D1]),  valable pour tous les $d$-tissus en   courbes dans une vari\'et\'e de dimension $n$. 

Ces tissus jouent un r\^ole important dans l'\'etude des $(n+3)$-tissus en courbes $W_{0,n+3}$  dits {\it exceptionnels},  qui   g\'en\'eralisent le 5-tissu de Bol ([Bo]) de la dimension 2 :  les $(n+1)$-sous-tissus de $W_{0,n+3}$ sont en effet tous de rang un.   Damiano, qui appelle    {\it combinatoires} les relations ab\'eliennes de $W_{0,n+3}$ que ces sous-tissus  engendrent, a d\'emontr\'e, du moins   si $n$ \'etait pair,  que  c'\'etait  tout l'espace 
des relations ab\'eliennes de $W_{0,n+3}$ qui \'etait ainsi obtenu, y compris  la relation ab\'elienne dite {\it d'Euler} d\'efinie par Gelfand et  MacPherson ([GM])  dans le contexte analytique r\'eel ([D1][D2]). Si $n$ est impair, Pirio  a d\'emontr\'e  qu'il n'en \'etait plus de m\^eme, et que l'espace des relations ab\'eliennes de $W_{0,n+3}$ \'etait alors somme directe de l'espace des relations ab\'eliennes combinatoires,  
et de l'espace engendr\'e par la relation d'Euler ([P]). C'est \`a cette occasion qu'il a pos\'e quelques questions, auxquelles nous r\'epondons partiellement dans cet article.
 
L'application $$Tr:(V_1,\cdots,V_{n+1})\to V_1+\cdots +V_{n+1}$$ de
$\bigoplus_{\lambda=1}^{n+1} T( {\cal F}_\lambda)$ dans $ T{\cal U}$ 
est de rang maximum $n$ en chaque point de $\cal U$, et  son noyau 
$$E:=Ker\Bigl(Tr: \bigoplus_{\lambda=1}^{n+1} T({\cal F}_\lambda)\to T{\cal U}\Bigr)$$est donc un fibr\'e vectoriel de rang 1. 

On se propose d'abord de d\'efinir  sur $\bigotimes^{n-1}E^*$ (la $(n-1)$-\`eme puissance tensorielle  du fibr\'e  dual \footnote{On peut aussi d\'efinir la connexion sur $\bigotimes^{n-1}E$ par dualit\'e.} de $E$), une connexion canonique (holomorphe ou analytique r\'eelle selon le contexte),     dont
 \pagebreak la \emph{courbure}   $\Omega$ g\'en\'eralise  pour  $n$ quelconque   la courbure de Blaschke ([B]) du cas $n=2$, le tissu \'etant de rang 1 ou 0 selon que $\Omega$ est nulle ou non.

On d\'ecrira ensuite, si $n=3$,  un proc\'ed\'e pour obtenir \emph{tous} les (germes de) 4-tissus en courbes de rang un, \`a isomorphisme pr\`es, \`a partir de trois fonctions $Q,u$ et $v$ sur la vari\'et\'e ambiante. Sans \^etre capable de pr\'eciser l'ensemble des triplets $(Q,u;v)$ engendrant des tissus isomorphes, 
nous montrerons cependant qu'il existe une infinit\'e de telles classes d'isomorphisme, contrairement \`a ce qui se paesse quand $n=2$ o\`u tous les 3-tissus de rang un sont (localement) isomorphes \`a un tissu form\'e de trois pinceaux de droites (relativement \`a la structure affine d\'efinie par des coordonn\'ees locales).

On appelle plus g\'en\'eralement \emph{quadrilat\'eral} \footnote{Il existe une d\'efinition s'appliquant \`a des situations beaucoup plus g\'en\'erales ([D1][D2]).}  un $(n+1)$-tissu en courbes dans un espace de  dimension $n$, localement isomorphe \`a $(n+1)$ pinceaux de droites dans un espace projectif de dimension $n$. Un tel tissu est de rang un. Mais,  pour $n\geq 3$, et contrairement \`a ce qui se passe pour $n=2$,  la quadrilat\'eralit\'e est une propri\'et\'e beaucoup plus forte que  le simple fait d'avoir un rang \'egal \`a un :  il se peut, par exemple,  qu'il existe des feuilletages en surfaces contenant deux des feuilletages ${\cal F}_\lambda$ et ${\cal F}_\mu$  du tissu,  et le plus grand nombre possible  $d$ de tels couples $(\lambda,\mu)$ est un invariant de $W$, g\'en\'eriquement \'egal \`a 0 pour un $(n+1)$-tissu en courbes, m\^eme s'il est de rang un, et \'egal \`a $\frac{n(n+1)}{2}$ pour un tissu quadrilat\'eral, avec tous les interm\'ediaires possibles\footnote{Pour $n\geq 4$, on pourrait aussi d\'efinir le plus grand nombre possible  de triplets $(\lambda,\mu,\nu)$  tels qu'il passe un feuilletage de codimension $n-3$ par ${\cal F}_\lambda$, ${\cal F}_\mu$ et ${\cal F}_\nu$, etc...}. 

\section{Courbure des $(n+1)$-tissus de courbes  en dimension $n$ }

\noindent Soit alors $\sigma=(V_\lambda)_\lambda $ une section locale partout non nulle  de $E$, $\bigl(\sum_{\lambda=1}^{n+1} V_\lambda=0\bigr)$, 
et
$\widehat{\sigma}$ la section duale de $E^*$ (d\'efinie comme la forme lin\'eaire sur le module des sections de  $E$ qui prend la valeur 1 sur $\sigma$).
On conviendra que

- les indices not\'es $\lambda,\mu,\nu...$ varient entre 1 et $n+1$, 

- ceux d'entre eux not\'es $i,j,k, ...$   varient entre 1 et $n$. 

\noindent La famille  $ (V_i)_{1\leq i\leq n} $ constitue alors une base locale du module des champs de vecteurs sur $\cal U$. On notera   
$ (\alpha_i)_{1\leq i\leq n} $ la base  locale duale du module des 1-formes (holomorphes ou analytiques r\'eelles suivant le cas), et
 $$V_{n+1}=-\sum_i V_i.$$
Si $[V_i,V_j]=\sum_k C_{ij}^k V_k$, $$d\alpha_k=-\sum_{i<j} C_{ij}^k\ \alpha_i\wedge\alpha_j.$$

\noindent Notons :

$A_i$ la $(n-1)$-forme\footnote{Il est sous-entendu que les indices $j$ dans les produits ext\'erieurs $\bigwedge_j$ ou $\bigwedge_{j,j\neq i}$ sont pris dans l'ordre croissant.} $$ A_i:=(-1)^{i+1}\ .\bigwedge_{j,j\neq i} \alpha_j $$

$\varphi_i$ la fonction
$$\varphi_i:=\sum_{j,j<i}C_{ji}^j-\sum_{j,i<j}C_{ij}^j$$

\noindent telle que 
$$dA_i=\varphi_i\ . \bigwedge_{j} \alpha_j,$$

  $\omega$ la 1-forme 
$$\omega:=\sum _i \varphi_i\ \alpha_i,$$de sorte que $$dA_i=A_i\wedge \omega,$$ 

\noindent et posons :$$\Omega:=d\omega.$$

\noindent {\bf Th\'eor\`eme 1 :}

{\it 

\noindent $(i)$ Les relations ab\'eliennes locales non nulles du tissu s'identifient naturellement aux fonctions $f$ telles que $$(*)\hskip .5cm d\ (log\ f)=-\omega .$$
Le tissu est donc de rang 1 ou 0 selon que la  2-forme  $\Omega$  est nulle ou non.

\noindent $(ii)$ La connexion $\nabla$ d\'efinie localement sur $\bigotimes^{n-1}E^*$ en posant  $$\nabla\bigl( \widehat{\sigma}^{\otimes (n-1)}\bigr) =   \omega .\bigl( \widehat{\sigma}^{\otimes (n-1)}\bigr)$$ne d\'epend ni de l'ordre des champs de vecteurs $V_\lambda$, ni de la section locale $\sigma=\{V_\lambda\}_\lambda$ de $E$ : les connexions ainsi d\'efinies localement se recollent et permettent de d\'efinir globalement  une connexion sur $\bigotimes^{n-1}E^*$.

Les relations ab\'eliennes du tissu s'identifient donc encore naturellement aux sections de $\bigotimes^{n-1}E^*$ \`a d\'eriv\'ee covariante nulle 
et le tissu est de rang 1 ou 0 selon que la courbure $\Omega$ de cette connexion   est nulle ou non.

\noindent $(iii)$ Pour $n=2$, la  courbure $\Omega$   est la courbure de Blaschke.
 }
  \vspace{3 mm}

 \noindent {\it D\'emonstration :}

Pour qu'une $(n-1)$-forme $\sum_j f_j.A_j$ soit ${\cal F}_i$-basique, il faut et il sufit que soient r\'ealis\'ees les deux conditions suivantes :

- les fonctions $f_j$ sont identiquement nulles si $j\neq i$, 

- la fonction $f_i$ est telle que $L_{V_i}(f_i.A_i)=0$, ce qui s'\'ecrit encore : $$(*_i)\hskip .5cm (V_i.f_i) +f_i\ \varphi_i=0.$$
\indent Pour que $\Bigl((f_1.A_1),(f_2.A_2),\cdots,(f_1.A_1),-\sum_i(f_i.A_i)\Bigr)$ soit une relation ab\'elienne, il faut et il sufit que soient r\'ealis\'ees les  conditions suivantes :

- toutes les fonctions $f_i$ sont   \'egales \`a une m\^eme fonction $f$, qui doit donc v\'erifier chacune des \'equations $(*_i)$, ce qui s'\'ecrit, si $f$ n'est pas nulle : $$(*)\hskip .5cm d\ (log\ f)=-\omega .$$
 
- la condition 
$L_{V_{n+1}}\eta_{n+1}=0$ doit \^etre r\'ealis\'ee, o\`u l'on a pos\'e : $\eta_{n+1}=\sum_i(f_i.A_i)$.

\noindent La premi\`ere condition est la transcription de la condition $\iota_{V_{n+1}}(\eta_{n+1})=0$.  
Quant \`a la seconde, qui s'\'ecrit encore :
$$\sum_{i,j} L_{V_{i}} (f.A_j)=0,$$
 elle est automatiquement v\'erifi\'ee d\`es lors que $f$ v\'erifie $(*)$ : en effet, 
 $$L_{V_{i}}(f\ A_j)=(V_i. f)\ A_j+f\Bigl(d(\iota_{V_i}A_j)+\varphi_j\ A_j\Bigr)$$
et  $\iota_{V_i}(A_j)=-\iota_{V_j}(A_i)$ pour $i\neq j$, tandis que  $\iota_{V_i}(A_i)=0$,  
d'o\`u la partie $(i)$ du th\'eor\`eme.

Calculons la forme $\omega':=\sum_j\psi_j\ \beta_j$ que l'on aurait obtenue \`a la place de $\omega=\sum_i\varphi_i\ \alpha_i$ si l'on avait modifi\'e l'ordre des indices $\lambda$ en permutant    $\lambda_o$  et $n+1$. On obtiendrait :

 $\beta_i=\alpha_i-\alpha_{\lambda_o}$ si $i<\lambda_o$ , 

 $\beta_{\lambda_o}=-\alpha_{\lambda_o}$ , 

  $\beta_i=\alpha_{i+1}-\alpha_{\lambda_o}$ si $\lambda_o<i\leq n$.

 \noindent On en d\'eduit : $$\bigwedge_{i} \beta_i=- \bigwedge_{i}\alpha_i$$
 tandis que 
 
$A_i$ devient $-A_i$, d'o\`u $\psi_i=-\varphi_i$  si $i\neq\lambda_o$,  

 $A_{\lambda_o}$ devient $-\sum_{i\neq \lambda_o} A_{i}$, 
 
 \noindent d'o\`u $\psi_{\lambda_o}=\sum_i\varphi_i$, et $\omega'=\omega$ : ceci prouve que la connexion d\'efinie localement sur $\bigotimes^{n-1}E^*$ en posant $$\nabla \bigl( \widehat{\sigma}^{\otimes {n-1}}\bigr)=    \omega .\bigl( \widehat{\sigma}^{\otimes {n-1}}\bigr)$$ne d\`'epend pas de l'ordre des champs de vecteurs d\'efinissant $\sigma$.
 
 Elle ne d\'epend pas non plus de la section $\sigma$ de $E$. En effet, si
 $u$ d\'esigne une fonction arbitraire partout non nulle et si l'on remplace $\sigma=\{V_\lambda\}_\lambda$ par
  $\sigma'=\{uV_\lambda\}_\lambda$, la base duale $(\alpha'_i)_i$ v\'erifie : $\alpha'_i=\frac{1}{u}\alpha_i$ et
  $A'_i\bigl(=\bigwedge_{j,j\neq i}\alpha'_j\bigr)=\frac{1}{u^{n-1}}A_i$. On en d\'eduit :
  $$dA'_i=\omega'\wedge\Bigl(\bigwedge_{j\neq i}\alpha'_j\Bigr) $$
  avec $$ \omega'= \omega + d\ log\Bigl(\frac{1}{ u^{n-1} }\Bigr) \ .$$
  Puisque $u$ est une fonction de transition du fibr\'e $E$ de rang un, $\frac{1}{u^{n-1}}$ est la fonction de transition  du fibr\'e $\bigotimes^{n-1}E^*$ qui lui correspond, 
  d'o\`u la partie $(ii)$ du th\'eor\`eme.
  
  Rappelons que, que si  $\eta$ d\'esigne une 1-forme int\'egrable, l'\'equation $d\eta=\eta\wedge \omega$ exprime que $\omega$ est la forme d'une connexion de Bott sur le fibr\'e normal $N_{\cal F}$ au feuilletage   $ \cal F$ de codimension  un d\'efini par $T({\cal F})=ker\ \eta$  (forme de connexion relative \`a la trivialisation de $N_{\cal F}$ d\'efinie par la projection  parall\`element \`a $T({\cal F})$  de  tout champ de vecteur $Y$ tel que $\eta(Y) =1$).
  Par cons\'equent, si $n=2$, les \'equations  
   $d\alpha_1=\alpha_1\wedge \omega$ et  $d\alpha_2=\alpha_2\wedge \omega$   expriment   que $\omega$ est la forme d'une connexion   sur le fibr\'e $T_3$ engendr\'e par $V_3$, qui est une connexion de Bott \`a la fois pour le  feuilletage engendr\'e par $V_2$ et pour  celui engendr\'e par $V_1$,    quand on identifie  $T_3$ \`a leur fibr\'e normal commun, et qu'on  trivialise celui-ci par   $V_3$. Or il n'existe qu'une seule connexion ayant cette propri\'et\'e, et sa courbure est pr\'ecis\'ement  la courbure de Blaschke, d'o\`u la partie $(iii)$.

  \rightline{QED}
  
 On appellera \emph{courbure du tissu} la courbure de la connexion $\nabla$, c'est-\`a-dire 
  la 2-forme $$ \Omega= d\omega$$qui a une signification intrins\`eque, ind\'ependante des choix de $\sigma=\{V_\lambda\}_\lambda$ et de l'ordre des champs $V_\lambda$. Les $(n+1)$ tissus \`a courbure nulle (c'est-\`a-dire de rang 1) seront encore dits \emph{plats}. 
   
\section{Description  des 4-tissus en courbes de rang un, en dimension trois} 

Notons  $(x,y,z)$ les coordonn\'ees locales et donnons nous trois fonctions $Q,u,v$ de $(x,y,z)$.  Posons :
$$f=Q'_x-u.v'_x \ ,\ \ g=Q'_y-u.v'_y \ \hbox { et }\ \ h=Q'_z-u.v'_z\ ,$$et supposons que trois des quatre   2-formes $df\wedge dx$, $dg\wedge dy$,  $dh\wedge dz$ et
$du\wedge dv$ sont  toujours lin\'eairement ind\'ependantes. Notons $W$ le 4-tissu en courbes form\'e par les quatre feuilletages admettant respectivement comme syst\`eme de fonctions basiques $$(f,x),\ (g,y), (h,z)\ \hbox { et }(u,v)  $$
\noindent {\bf Th\'eor\`eme 2 :}{\it

\noindent   $(i)$ Le 4-tissu $W$ est plat.

\noindent   $(ii)$ R\'eciproquement,  tout germe de 4-tissu en courbes, dans une vari\'et\'e de dimension trois,  est n\'ecessairement isomorphe, s'il   est plat, \`a un tissu fabriqu\'e de cette fa\c con \`a partir de trois fonctions $Q,u,v$.} 
 \vspace{2 mm}

 \noindent {\it D\'emonstration : }
 
De la d\'efinition des fonctions $f,g,h$, on d\'eduit  les \'egalit\'es $$dQ=f.dx +g.dy+h.dz+u.dv $$et
 $$df\wedge dx+dg\wedge dy+dh\wedge dz+du\wedge dv=0.$$
 Comme les 2-formes $df\wedge dx$, $dg\wedge dy$, $dh\wedge dz$ et $du\wedge dv$ sont respectivement basiques relativement \`a  chacun des feuilletages constituant le tissu $W$, la deuxi\`eme de ces \'egalit\'es d\'efinit une 2-relation ab\'elienne non nulle sur le tissu, qui est donc de rang un.
 
 R\'eciproquement, on se donne  un germe de  4-tissu de codimension deux au voisinage d'un point en dimension trois. On peut toujours choisir les coordonn\'ees  locales  $(x,y,z)$ de fa\c con que $x$  (resp. $y$,   resp. $z$) soit fonction basique du premier  (resp. du second, resp. du troisi\`eme) feuilletage : il existe donc cinq fonctions $f,g,h,u,v$ de $(x,y,z)$, telles que les quatre feuilletages soient respectivement d\'efinis localement par les couples de fonctions basiques $(f,x)$, $(g,y)$, $(h,z)$ et $(u,v)$. 
 Supposons qu'il existe une relation ab\'elienne non triviale 
 $$A(f,x)\ df\wedge dx+B(g,y)\ dg\wedge dy+C(h,z)\ dh\wedge dz+D(u,v)\ du\wedge dv=0.$$
Notons  $$F(r,x)=\int A(r,x)\ dr  $$ $$G(r,y)=\int B(r,y)\ dr$$ $$  H(r,z)=\int C(r,z)\ dr$$ $$ U(r,v)=\int D(r,v)\ dr$$  une primitive de $A(r,x)$ $\Bigl($resp. de  $B(r,y)$, resp. de  $C(r,z)$, resp. de  $D(r,v)\Bigr)$ relativement \`a la variable $r$ :  on peut remplacer $f(x,y,z)$ par $F\bigl(f(x,y,z),x\bigr)$ dans la d\'efinition du premier feuilletage, et de m\^eme prendre $G\bigl(g,y\bigr)$ 
 au lieu de $g$,   $H\bigl(h,z\bigr)$  au lieu de $h$, et $U (u,v )$ au lieu de $u$. 
 
 \noindent Autrement dit, on peut toujours choisir les fonctions $f,g,h$ et $u$ de fa\c con que la relation ab\'elienne s'\'ecrive :
$$df\wedge dx+dg\wedge dy+ dh\wedge dz+ du\wedge dv=0, $$ce qui implique :
$$d\Bigl(f\ dx+g\ dy+h\ dz+u\ dv\Bigr)=0.$$
Au niveau des germes au voisinage d'un point, on peut appliquer le lemme de Poincar\'e : il existe une fonction $Q$ telle que 
$$dQ=f\ dx+g\ dy+h\ dz+u\ dv,$$
d'o\`u $$f=Q'_x-u.v'_x,\ g=Q'_y-u.v'_y,\ h=Q'_z-u.v'_z. $$

\rightline{QED}
 
Notons  que des triplets  $(Q_1,u_1,v_1)$ et $(Q_2,u_2,v_2)$ distincts peuvent d\'efinir des tissus isomorphes. Nous ne savons pas trouver une forme normale pour $(Q,u,v)$ caract\'erisant enti\`erement la classe d'isomorphisme du tissu associ\'e. 

Cependant, au voisinage de l'origine,  on peut faire un changement de coordonn\'ees de fa\c con  que les champs de vecteurs $(V_\lambda)$, $(1\leq \lambda\leq 4)$ d\'efinissant les feuilletages soient respectivement \'egaux \`a l'origine aux vecteurs $V_1^o=(0,0,1)$, $V_2^o=(1,0,0)$, $V_3^o=(0,1,0)$, $V_4^o=(1,1,1)$. On peut aussi imposer aux fonctions   $f, g, h, u, v$, toutes basiques pour l'un des   feuilletages du tissu, d'\^etre nulles \`a l'origine. On v\'erifie que ces conditions impliquent :

\noindent {\bf Lemme}: {\it 

On peut supposer, sans perte de g\'en\'eralit\'e \`a isomorphisme pr\`es :
$$ Q= \frac{x^2}{2}+yz +...,\hskip 1cm \ \ u=y-x+..., \ \ \hskip 1cm  v=z-x+...$$
les points de suspension indiquant des termes d'ordre sup\'erieur. 
}
 \vspace{2 mm}

 Le cas $\bigl(Q=yz+\frac{x^2}{2}$, $u=y-x$, $v=z-x\bigr)$ correspond au tissu dit \emph{quadrilat\'eral}, consistant en quatre pinceaux de droites parall\`eles\footnote{Ce tissu est aussi localement isomorphe (non projectivement) au tissu d\'efini par les quatre pinceaux de droites de sommets respectifs $[1:0:0:0]$, $[0:1:0:0]$, $[0:0:1:0]$, et $[0:0:0:1]$.} dans l'espace projectif de dimension 3. 
 \vspace{1 mm}

\noindent {\bf Un premier invariant (entier $d$ compris entre 0 et 6) :} 
 \vspace{1 mm}

On d\'efinit $d$ comme le  nombre   de couples   $1\leq \lambda < \mu\leq 4$ tels que 
$$V_\lambda\wedge V_\mu\wedge [V_\lambda,V_\mu]\equiv 0,$$c'est-\`a-dire tels qu'il existe un feuilletage ${\cal F}_{\lambda\mu}$ en surfaces contenant ${\cal F}_{\lambda }$ et ${\cal F}_{\mu}$.
Il est clair que cet entier est un invariant\footnote{Il est clair que des invariants analogues peuvent se d\'efinir dans des situations beaucoup plus g\'en\'erales.} de la classe d'isomorphisme du tissu.

\noindent {\it Exemple :} Notons $W_{P}$ le tissu, d\'ependant d'un param\^etre $P=(a,b,c,aa,bb,cc,e)$ form\'e de sept scalaires, d\'efini par :
$$u=y-x,\ v=z-x,$$ $$ Q=yz+\frac{x^2}{2}+a\ \frac{x^2y}{2} +aa\ \frac{y^2x}{2}+b\ \frac{y^2z}{2}+bb\ \frac{z^2y}{2}+c\ \frac{z^2x}{2}+cc\ \frac{x^2a}{2}     +e\ xyz.$$
Calculant $f:=Q'_x-uv'x$, $g:=Q'_y-uv'y$,  et $h:=Q'_z-uv'z$, 
et notant en abr\'eg\'e $\partial _x, \partial _y$ et $ \partial _z$ les champs de vecteurs coordonn\`ees, on peut alors prendre :
$$\hskip -5.6 cm V_1=f'_z\ \partial _y-f'_y\ \partial _z,$$
$$\hskip -5.3cm V_2=-g'_z\ \partial _x+ g'_x\ \partial _z,$$
$$\hskip -5.5cm V_3=h'_y\ \partial _x-h'_x\ \partial _y,$$
$$\hskip -5.7cm V_4=   \partial _x + \partial _y+  \partial _z.$$
Notant $D_{\lambda\mu}$ le d\'eterminant $V_\lambda\wedge V_\mu\wedge [V_\lambda,V_\mu]$, on   v\'erifie alors :
  
 $d=0$ \ \ dans le cas {\it g\'en\'erique}, 
 
 $d=1$ \ \  pour $P=( 0,0,-2,-1,0,1,1)$ :\hskip 1cm seul, $D_{14}$ est nul,
 
  $d=2$ \ \  pour $P=( 0,0,0,1,1,0,0)$ : \hskip 1.4 cm seuls, $D_{12}$ et $D_{14}$ sont nuls,
  
  $d=3$ \ \  pour $P=( -2,-2,-2,1,1,1,1)$ :\hskip .7cm seuls, $D_{14}$,  $D_{24}$ et $D_{14}$ sont nuls,
  
  $d=4$ \ \  pour $P=( -1,-1,0,1,1,0,0)$ :\hskip 1cm seuls, $D_{12}$, $D_{14}$,  $D_{24}$ et $D_{14}$ sont nuls,

 $d=5$ \ \  pour $P=( -1,0,0,1,0,0,0)$ : \hskip 1.2cm seul, $D_{23}$ n'est pas nul,
 
  $d=6$ \ \  pour $P=( 0,0,0,0,0,0,0)$ :\hskip 1.5cm c'est le cas {\it quadrilat\'eral}.

\noindent Ceci prouve d\'ej\`a (contrairement \`a ce qui se passe en dimension deux o\`u tous les 3-tissus de rang un   sont localement isomophes) :  

\indent -  que la quadrilat\'eralit\'e est une propri\'et\'e beaucoup plus forte, pour un tel 4-tissu,  que le simple fait d'\^etre de rang un,   \indent

 - et qu'il existe plusieurs classes d'isomorphisme  dans l'ensemble des  4-tissus de courbes  qui sont de rang un, en dimension trois ; il  en existe en fait une infinit\'e :

\noindent {\bf Th\'eor\`eme 3 : }

 {\it Il existe une infinit\'e non d\'enombrable de classes d'isomorphisme de germes de 4-tissus en courbes, de rang un,  dans une vari\'et\'e de dimension trois.}
\vspace{1 mm}

 \noindent C'est un corollaire imm\'ediat de la proposition suivante, dans laquelle 
 on note $\phi_{ijk}$ les coefficients de la s\'erie de Taylor de toute  fonction $\phi$ au voisinage de l'origine, 

\noindent {\bf Proposition :}

{\it D\`es lors que son d\'enominateur n'est pas nul,  l'expression } $$\frac{2Q_{120}-v_{110}+2v_{020}}{2Q_{201}+u_{101}+v_{101}}$$
{\it est  un invariant de la classe d'isomorphisme du tissu, que l'on peut faire varier \`a volont\'e.}
 \vspace{1 mm}

\noindent {\it D\'emonstration :}

Supposons qu'il existe un changement de coordonn\'ees $(x,y,z)\mapsto (X,Y,Z)$ r\'ealisant un  isomorphisme de $W$ (d\'efini par $(Q,u,v)$ sur
$W'$ (d\'efini par $(QQ,uu,vv))$, avec 
$$\hskip -1cm u:=y-x+...,\hskip  .3cm\ \ v:=z-x+..., \hskip  .5cm Q:=yz+\frac{x^2}{2}+...,$$et $$uu:=Y-X+...,\ vv:=Z-X+..., \ \ QQ:=YZ+\frac{X^2}{2}+.... .$$
On en d\'eduit :

\noindent $f=Q'_x-uv'_x$ et $F=QQ'_X-UV'_X$,

\noindent   $g=Q'_y-uv'_y$,     et $G=QQ'_Y-UV'_Y$, 

\noindent $h=Q'_z-uv'_z$ et $H=QQ'_Z-UV'_Z$.

\noindent Puisque $X$ et $F$ doivent \^etre des fonctions basiques du   feuilletage image de celui d\'efini par   $(x,f)$, il doit exister des fonctions de deux  variables $R$ et $M$ telles que 

\noindent $X=R(f,x)$ et $F(X,Y,Z)=M(f,x)$.

\noindent Il existe de m\^eme des fonctions de deux variables $S,T,N,P, U,V$ telles que

\noindent $Y=S(g,y)$ et $G(X,Y,Z)=N(g,y)$,

\noindent $Z=T(h,z)$ et $H(X,Y,Z)=P(h,z)$,

\noindent  $uu(X,Y,Z)=U(u,v)$ et $vv(X,Y,Z)=V(u,v)$.

On obtient donc cinq identit\'es :
$$F\Bigl(R\bigl(f(x,y,z),x\bigr),S\bigl(g(x,y,z),y\bigr),T\bigl(h(x,y,z),z\bigr)\Bigr)\equiv M\bigl(f(x,y,z),x\bigr),$$
$$G\Bigl(R\bigl(f(x,y,z),x\bigr),S\bigl(g(x,y,z),y\bigr),T\bigl(h(x,y,z),z\bigr)\Bigr)\equiv N\bigl(g(x,y,z),y\bigr),$$
$$H\Bigl(R\bigl(f(x,y,z),x\bigr),S\bigl(g(x,y,z),y\bigr),T\bigl(h(x,y,z),z\bigr)\Bigr)\equiv P\bigl(h(x,y,z),z\bigr),$$
$$ uu\Bigl(R\bigl(f(x,y,z),x\bigr),S\bigl(g(x,y,z),y\bigr),T\bigl(h(x,y,z),z\bigr)\Bigr)\equiv U\bigl(u(x,y,z),v(x,y,z)\bigr),$$
$$ vv\Bigl(R\bigl(f(x,y,z),x\bigr),S\bigl(g(x,y,z),y\bigr),T\bigl(h(x,y,z),z\bigr)\Bigr)\equiv V\bigl(u(x,y,z),v(x,y,z)\bigr).$$
D\'eveloppant ces identit\'es en s\'eries de Taylor en $(x,y,z)$ au voisinage de l'origine,  
on constate  d'abord qu'\`a l'ordre 1,   les contraintes impos\'ees \`a $(Q,u,v)$ et $(QQ,uu,vv)$ impliquent l'existence d'un scalaire $m$ tel que  les fonctions  $R,S,T,M,N,P, U,V$ de deux variables $(r,s)$ aient, \`a l'origine,  des jacobiens  tous \'egaux \`a $m.Id_2$ :

On obtient alors, \`a l'ordre 2 et \`a l'aide de Maple,  un certain nombre de relations, et en particulier les deux relations $$2Q_{201}+u_{101}+v_{101}=m(2QQ_{201}+uu_{101}+vv_{101}) $$ et $$2Q_{120}-v_{110}+2v_{020}=m(2QQ_{120}-vv_{110}+2vv_{020}) .$$
 Ceci prouve que $$\frac{2Q_{120}-v_{110}+2v_{020}}{2Q_{201}+u_{101}+v_{101}}$$
est  un invariant de la classe d'isomorphisme du tissu, que l'on peut faire varier \`a volont\'e d\`es lors que le d\'enominateur n'est pas nul.

\rightline{QED}

\noindent {\bf Ref\'erences}

\noindent [B] W.Blaschke : {\it \"Uber Geweben von Kurven in $\reals^3$}, Abh. Math.Hamburg 11, 1936, 387-393.

\noindent [Bo] G.Bol : {\it \"Uber ein bemerkenswertes F\"unfgewebe in der Ebene}, Abh. Math. Sem. Univ. Hamburg 9, 1933, 291-298.

\noindent [D1] D.Damiano : {\it Webs,abelian equations and characteristic classes},PhD thesis Brown Univ., 1980.

\noindent [D2] D.Damiano : {\it Webs  and characteristic forms on Grassmann manifolds},Am.J. of Maths.105, 1983,  1325-1345.

\noindent [GM] I.Gelfand and R. MacPherson : {\it Geometry in Grassmannians and a generalization of the  dilogarithm}, Adv. in Math. 44, 1982, 279-312.

\noindent [P] L. Pirio : {\it On the (n+3)-webs  by rational curves induced by the forgetful maps on the moduli spaces ${\cal
 M}_{0,n+3}$}, arXiv 2204.04772.v1, [Math AG], 10-04-2022.
 
 \end{document}